\documentclass{cedram-aif}
\usepackage{amsmath,amsfonts,amsthm,amssymb,color, enumerate}
\numberwithin{equation}{section}

\newtheorem{corollary}[equation]{Corollary}

\theoremstyle{remark}
\newtheorem{remark}[equation]{Remark}
\theoremstyle{definition}
\newtheorem{defn}[equation]{Definition}

\newcommand{\tr}{\mathop{\rm tr}}
\newcommand{\cD}{{\mathcal{D}}}
\newcommand{\cN}{{\mathcal{N}}}
\newcommand{\cP}{{\mathcal{P}}}
\newcommand{\cR}{{\mathcal{R}}}
\newcommand{\N}{\mathbb{N}}
\newcommand{\bbH}{\mathbb{H}}
\newcommand{\cC}{\mathcal{C}}
\newcommand{\cL}{\mathcal{L}}
\newcommand{\Q}{\mathbb{Q}}
\newcommand{\fp}{\mathfrak{p}}
\newcommand{\Vol}{\textrm{Vol}}
\newcommand{\pa}{\partial}
\newcommand{\rint}{\sideset{^{0}}{}\int}
\newcommand{\res}{\textrm{res }} 
\newcommand{\Res}{\textrm{Res }} 
\newcommand{\Tr}{\textrm{Tr }}

\def\BIgl({\mathopen{\hbox{\smaller${\biggl(}$}}}
\def\BIgr){\mathclose{\hbox{\smaller${\biggr)}$}}}
\def\BIglv|{\mathopen{\hbox{\smaller${\biggl|}$}}}
\def\BIgrv|{\mathclose{\hbox{\smaller${\biggr|}$}}}
\def\bil({\mathopen{\raise.25pt\hbox{\smaller${\bigl(}$}}}
\def\bir){\mathclose{\raise.25pt\hbox{\smaller${\bigr)}$}}}
\def\bilv|{\mathopen{\raise.25pt\hbox{\smaller${\bigl|}$}}}
\def\birv|{\mathclose{\raise.25pt\hbox{\smaller${\bigr|}$}}}

\renewcommand{\Re}{\mathop{\rm Re}}

\def\udstrutpt#1#2{\noindent\vbox to #1pt{}\noindent\lower#2pt\vbox{}\ignorespaces}

\def\zeroint_#1{\mathchoice
{\sideset{ ^{\hskip5pt 0\hskip-5pt}}{}\int\nolimits_{\hskip-5pt #1}}
{^{\hskip4pt 0\hskip-4pt}\int_{#1}} 
{\sideset{ ^{\hskip5pt 0\hskip-5pt}}{}\int}
{\sideset{ ^{\hskip5pt 0\hskip-5pt}}{}\int}
}

\def\smint{\raise.5pt\hbox{\smaller$\dsty \int$}}
{\makeatletter
\gdef\citemt[#1]#2{[\hyper@@link[cite]{}{cite.#2}{#1}]}
\gdef\citelcs[#1]#2{(loc.\ cit.)}
\gdef\citelc[#1]#2{(\hyper@@link[cite]{}{cite.#2}{loc.\ cit.})}
\gdef\citeib[#1]#2{[\hyper@@link[cite]{}{cite.#2}{ibidem}, #1]}
\gdef\citeibna#1{(\hyper@@link[cite]{}{cite.#1}{ibidem})}
}

\def\udstrutpt#1#2{\noindent\vbox to #1pt{}\noindent\lower#2pt\vbox{}\ignorespaces}

\renewcommand{\qed}{\hspace*{\fill} \setlength{\unitlength}{1mm}
\begin{picture}(2.5,2.5)
      \put(0,0){\framebox(2.5,2.5){}}
  \end{picture}
\setlength{\unitlength}{1pt}}

\newcommand{\dvol}{\, \textrm{dvol}}

\newcommand{\bbS}{\mathbb{S}}

\newcommand{\R}{\mathbb{R}}
\newcommand{\C}{\mathbb{C}}

\newcommand{\cS}{\mathcal{S}}

\title[On the spectral theory and dynamics of a.h. manifolds]{On the spectral theory and dynamics of asymptotically hyperbolic manifolds}

\alttitle{Sur la th\'eorie spectrale et la dynamique des vari\'et\'es asymptotiquement hyperboliques}

\author{\firstname{Julie} \lastname{Rowlett}}

\address{Hausdorff Center for Mathematics\\ Villa Maria
Endenicher Allee 62 \\ D-53115 Bonn, Deutschland} 

\email{rowlett@math.uni-bonn.de}

\keywords{asymptotically hyperbolic, conformally compact, wave trace, negative curvature, resonances, length spectrum, topological entropy, dynamics, geodesic flow, prime orbit theorem, quantum and classical mechanics.}
  
\altkeywords{vari\'et\'e asymptotiquement hyperbolique, conformement compact, courbures negatives, spectre des longeurs g\'eodesiques, flot g\'eodesique, dynamique, formule de trace dynamique, entropie topologique, m\'ecanique quantique et classique.}

\subjclass{37D40, 58J50, 53C22}

\begin{document}

\maketitle

\begin{abstract} %%%%%%%%%%%%%%%%%%%
We present a brief survey of the spectral theory and dynamics of infinite volume asymptotically hyperbolic manifolds.  Beginning with their geometry and examples, we proceed to  their spectral and scattering theories, dynamics, and the physical description of their quantum and classical mechanics.  We conclude with a discussion of recent results, ideas, and conjectures.  
\end{abstract}%%%%%%%%%%%%%%%%%%%

\begin{altabstract}
Cet article est une pr\'esentation rapide de la th\'eorie spectrale et de la dynamique des vari\'et\'es asymptotiquement hyperboliques \`a volume infini.   Nous commen\c cons  par leur g\'eom\'etrie et quelques exemples, on poursuit en rappelant  leur th\'eorie spectrale, puis continuons sur des d\'eveloppements r\'ecents de leur dynamique.  Nous concluons par une discussion de r\'esultats qui d\'emontrent un rapport entre leurs m\'ecaniques quantiques et classiques et enfin, nous offrons quelques id\'ees et conjectures.
\end{altabstract}

\section{Introduction and contents}

In 1984, Fefferman and Graham \cite{fg} introduced the concept of a conformally compact metric to generalize the Poincar\'e model of hyperbolic space and study conformal invariants.  These metrics are also known as Poincar\'e metrics.  A Riemannian manifold with boundary equipped with a conformally compact metric is known as a conformally compact manifold.  In 1987, Mazzeo and Melrose \cite{mmah} defined a special class of conformally compact manifolds, the \em asymptotically hyperbolic manifolds.  \em   These are compact manifolds with boundary equipped with an asymptotically hyperbolic metric:  a complete metric with sectional curvatures asymptotically equal to $-1$ at infinity.   Although the asymptotically hyperbolic manifolds are a proper subset of the conformally compact manifolds, they are a natural generalization of infinite volume hyperbolic manifolds and contain several interesting examples including manifolds with Poincar\'e metrics which are also Einstein.  These Poincar\'e Einstein metrics arise in AdS-CFT correspondence in string theory.     

Both finite and infinite volume hyperbolic manifolds have attracted interest from analytic number theory and arithmetic geometry, mathematical physics and spectral theory, and dynamics.  In many cases, these spaces admit a Selberg trace type formula \cite{sel}.  Such a trace formula demonstrates a relationship between their spectral theory and dynamics.  Phillips \cite{ph} and Patterson \cite{pat} obtained a precise relationship between the spectrum and the topological entropy of the geodesic flow.  In the language of mathematical physics, these results describe interactions between the quantum and classical mechanics.  

Motivated by the growing interest of mathematicians and physicists in conformally compact manifolds in general and Poincar\'e Einstein manifolds in particular, since all Poincar\'e Einstein metrics are asymptotically hyperbolic, it is natural to pursue a deeper understanding of the dynamics and its relationship to the spectral theory of these manifolds.  The purpose of this survey is:  to collect and summarize references for infinite volume asymptotically hyperbolic manifolds, to present the spectral theory and dynamics of these spaces, and to present recent developments concerning interactions between their quantum and classical mechanics.  This work is organized as follows.  

\begin{enumerate}
\item Introduction and contents%%%%%
\item Geometry %%%%%%
\begin{enumerate}[i.]
\item Model geometry
\item Asymptotic geometry
\item Renormalized integral
\end{enumerate}
\item Spectral and scattering theories %%%%%%%%%%%
\begin{enumerate}[i.]
\item The resolvent
\item The scattering operator
\item Counting estimates
\end{enumerate}
\item Dynamics %%%%%%%%%%%%
\begin{enumerate}[i.]
\item Dynamics of asymptotically hyperbolic manifolds
\item Dynamical zeta functions 
\end{enumerate} 
\item Quantum and classical mechanics %%%%%%%%%%%
\begin{enumerate}[i.]
\item The wave group
\item Trace formulae
\item Interactions between quantum and classical mechanics
\end{enumerate}
\item The horizon %%%%%%%%%%%%%%%
\begin{enumerate}[i.]
\item Wave trace remainder estimates
\item Spectral zeta functions 
\end{enumerate}
\end{enumerate} 

\subsection*{Bibliographical notes} 
Fefferman and Graham \cite{fg} defined conformally compact metrics to study conformal invariants.  A conformally compact manifold does not have a unique metric at the boundary; rather, the metric induces a conformal class of metrics at the boundary which is known as the ``conformal infinity.''  For a real analytic conformal Riemannian manifold, they showed that one may associate a conformally compact manifold in one dimension higher, such that the conformal manifold is the conformal infinity.  Based on this ambient metric construction, \cite{fg} computed several scalar curvature invariants for conformal Riemannian manifolds.  

Mazzeo and Melrose studied the resolvent operator for asymptotically hyperbolic metrics in \cite{mmah}; their main results will be reviewed in \S 3.1.  Mazzeo went on to study the $\cL^2$ cohomology of conformally compact manifolds and proved their Hodge theorem in \cite{ccma}.  

Selberg proved trace formulae for (weakly) symmetric spaces using a beautiful combination of techniques from harmonic analysis, group theory, and representation theory \cite{sel}.  The trace formula has been generalized to various other settings including congruence subgroups and $PSL(2,\R)$ by Hejhal \cite{he}, \cite{he2}; and some spaces of infinite volume by Gangolli and Warner \cite{gw}, \cite{gw2}.  Moreover, interactions between the trace formula and Hecke operators in arithmetic geometry have been investigated by Arthur \cite{art}.

Phillips \cite{ph} and Patterson \cite{pat} demonstrated that the existence of pure point spectrum determines and conversely is determined by the topological entropy of the geodesic flow; see also joint work of Phillips and Sarnak \cite{phs}.  

\section{Geometry}%%%%%%%%%%%
The model geometry of an asymptotically hyperbolic manifold is an infinite volume hyperbolic manifold which is the quotient of real hyperbolic space by a \em convex cocompact \em group. 

\subsection{Model geometry}
\begin{defn}  A finitely-generated discrete torsion-free group $\Gamma$ of isometries of real hyperbolic space is \em convex cocompact \em if $\Gamma$ does not contain parabolic elements, contains at least one hyperbolic element, and admits a finite sided fundamental domain with infinite volume.  
\end{defn}  

Examples of convex cocompact groups include Schottky groups, Fuchsian groups, and quasifuchsian groups; the term \em convex cocompact \em first appeared in Sullivan's work \cite{su}.  However, in that work, convex cocompact groups are defined to be those groups which do not contain parabolic elements and \em may \em contain hyperbolic elements.  In particular, he showed that the quotient of real hyperbolic space by a convex cocompact group is either compact or has infinite volume.  

\begin{defn} For a discrete torsion-free group $\Gamma$ of isometries of $n+1$ dimensional real hyperbolic space $\bbH^{n+1}$, the quotient $\bbH^{n+1} / \Gamma$ is \em convex cocompact \em if $\Gamma$ is a convex cocompact group.  The \em limit set \em of $\Gamma$, $\Lambda_{\Gamma}$, is the smallest closed $\Gamma$-invariant subset of $\pa \bbH^{n+1} = \bbS^n$; it is equivalently given by the intersection of the closure of the orbit by $\Gamma$ of any point of $\bbH^{n+1}$ with $\pa \bbH^{n+1}$.  The \em convex core \em is defined by 
$$CH(\Lambda_\Gamma)/ \Gamma \subset \bbH^{n+1} / \Gamma,$$
where $CH(\Lambda_{\Gamma})$ is the convex hull of $\Lambda_{\Gamma}$; for a convex cocompact group $\Gamma$, the convex core is a compact subset of $\bbH^{n+1} / \Gamma$.  
\end{defn}  

For a discrete, finitely-generated, torsion-free group $\Gamma$ of isometries of $\bbH^{n+1}$, the quotient is conformally compact precisely when $\Gamma$ is convex cocompact as defined above; see for example \cite{bor}.  We recall the definition of a conformally compact metric below.  

\subsection{Asymptotic geometry}

Asymptotic geometries such as asymptotically hyperbolic, asymptotically conic, and asymptotically cylindrical are Riemannian manifolds with boundary whose Riemannian metric admits a particular asymptotic form near the boundary; the definitions of these spaces are primarily due to Melrose \cite{tapsit}.  The definitions of all the aforementioned asymptotic geometries and the definition of conformally compact are formulated in terms of a \em boundary defining function.  \em 

\begin{defn}  Let $(X^{n+1}, \pa X)$ be a smooth, compact $n+1$ dimensional manifold with boundary.  A \em boundary defining function \em $x$ is a smooth function defined in a neighborhood $\cN$ of $\pa X$ so that $x: \cN \to [0, \infty)$, $\pa X = x^{-1} (\{0\}),$ and $dx \neq 0$ on $\pa X$.  \end{defn}

A conformally compact manifold is a Riemannian manifold with boundary whose metric admits the following form.  
\begin{defn} Let $X^{n+1}$ be a smooth compact manifold with boundary $\pa X$ and boundary defining function $x$.  If there exists a smooth metric $\bar{g}$ defined on $X$ which is non-degenerate up to $\pa X$ such that the metric $g$ on the interior of $X$ has the form 
$$g = \frac{\bar{g}}{x^2},$$
then $g$ is a conformally compact metric. 
\end{defn}

The following definition is from \cite{bor} and is related to the geometric setting of earlier work by Guillop\'e and Zworski \cite{gzpol}.  
\begin{defn} A conformally compact manifold $(X^{n+1}, g)$ is \em hyperbolic near infinity \em if there exists a (possibly disconnected) convex cocompact hyperbolic manifold $(X_0, g_0)$ and compact sets $K \subset X$, $K_0 \subset X_0$, such that 
$$(X-K, g) \cong (X_0 - K_0, g_0).$$
\end{defn} 
These spaces also appear in the literature as ``compact perturbations of convex cocompact hyperbolic manifolds.''  They enjoy many of the beautiful properties of their unperturbed model spaces; \cite{bor} proved a Poisson formula and estimates for their spectral counting function.  

The definition of an asymptotically hyperbolic manifold is originally due to Mazzeo and Melrose \cite{mmah}.  

\begin{defn}  A manifold with boundary $(X^{n+1}, \partial X)$ is \em asymptotically hyperbolic \em if there exists a boundary defining function $x$ so that in a non-empty neighborhood of the boundary $\cN \cong (0, x_0) \times \partial X,$ the metric 
\begin{equation}\label{ahmetric} g = \frac{dx^2 + h(x, y, dx, dy)}{x^2} \end{equation}
where $h|_{\{x=0\}}$ is independent of $dx$ and is a Riemannian metric on $\pa X$. 
\end{defn}  

\begin{remark} It is straightforward to demonstrate that every asymptotically hyperbolic manifold is conformally compact.  
%For an asymptotically hyperbolic manifold $(X, g)$, let $x$ be a boundary defining function such that on a neighborhood $\cN \cong \pa X \times (0, x_0)$ 
%$$g = \frac{dx^2 + h(x, y, dy)}{x^2} \textrm{ on $\cN$.}$$  
%Let $\chi \in \cC^{\infty} (X)$, $0 \leq \chi \leq 1$ satisfy 
%$$\chi = \left\{ \begin{array}{ll} 0 & \textrm{ on } \cN' \cong \left( 0, \frac{x_0}{2} \right) \times \pa X \\ 1 & \textrm{ on } X - \cN \end{array} \right .$$
%Let 
%$$\rho := (1 - \chi)x + \chi.$$
%Then, $\rho \equiv x$ on $\cN'$, so $\rho$ is also a boundary defining function. 
%$$g = \frac{d\rho^2 + h(\rho, y, dy)}{\rho^2} \qquad \textrm{on } \cN'.$$Moreover, $\rho \equiv 1$ on $X - \cN$, so $\rho^2 g \equiv g$ on $X - \cN$.%%
%By definition of $\rho$ and $x$, $0 < \rho$ on $X$.  Consequently, $\bar{g} := \rho^2 g$ is a Riemannian metric on $X$ which is non-degenerate up to $\pa X$, and 
%$$g = \frac{\bar{g}}{\rho^2}.$$
On the other hand, the sectional curvatures of conformally compact manifolds approach $-|dx|_{x^2 g} ^2$ at the boundary:  this is not necessarily a constant function.  Thus, not all conformally compact manifolds are asymptotically hyperbolic; see \cite{bor2}.  
\end{remark} 

Recall that a Riemannian metric is \em Einstein \em if the Ricci and metric tensors satisfy the relationship 
$$\textrm{Ric} = - c g,$$
for some constant $c$.  In the case of Poincar\'e Einstein metrics, those which are both conformally compact and Einstein, $c > 0$, and is usually normalized to $c=n+1$ in dimension $n+1$.  Some explicit examples of the Poincar\'e Einstein metrics which arise in AdS/CFT correspondence in string theory include the hyperbolic analogue of the Schwarzschild metric, and in dimension $n+1 = 4$, the Taub-BOLT metrics on disk bundles over $\bbS^2$.  These examples may be found in Anderson \cite{and2}.  

Mazzeo and Melrose observed that along any smooth curve in $X - \partial X$ approaching a point $p \in \partial X,$ the sectional curvatures of $g$ approach $- |d x|^2 _{x^2 g}.$  For each $h$ in the conformal class of $h|_{\pa X}$, by \cite{mmah} and \cite{mcpe}, there exists a unique (near the boundary) boundary defining function $x$ so that 
\begin{equation} \label{eq:seccurv} |dx|_{x^2 g} = 1, \quad \textrm{and} \quad g = \frac{dx^2 + h(x, y, dy)}{x^2} \textrm{ near $\partial X$.} \end{equation} 
With this normalization, the sectional curvatures are $-1+ O(x)$ as $x\to 0$, so it is natural to call these spaces ``asymptotically hyperbolic."  If the metric is also Einstein with Ric $g = - (n+1)g$, then the sectional curvatures are  $-1 + O(x^2)$ in a neighborhood $\pa X$.  When the sectional curvatures converge to $-1$ at a faster rate, one has rigidity results; see Shi and Tian \cite{shti}, and Anderson \cite{and}.  The rate at which the sectional curvatures converge to $-1$ was also used by Bahuaud in  \cite{be-th} and \cite{be} to give an intrinsic condition for Lipschitz conformally compact asymptotically hyperbolic manifolds.  

The reader is advised that the term ``asymptotically hyperbolic'' has been used differently by other authors including M\"uller \cite{mu}.  However, those ``asymptotically hyperbolic'' manifolds are not  conformally compact.  The asymptotically hyperbolic manifolds considered here are those which are also conformally compact, but for the exposition, we shall refer to them simply as asymptotically hyperbolic.  

\subsection{Renormalized integral}  
Since conformally compact manifolds have infinite volume, it is natural to introduce an integral renormalization.  Recall that the finite part $ \underset{{\epsilon=0}}{\mathrm{f.p.}} f(\epsilon)$ is defined as $f_0$ when 
$$f(\epsilon) = f_0 + \sum_{k=1} ^{\infty} f_k \epsilon^{- \lambda_k} (\log{\epsilon})^{m_k} + o(1),$$ 
with $\mathfrak{R}( \lambda_k ) \geq 0$, and $m_k \in \N \cup \{0\}.$  Uniqueness of $f_0$ is demonstrated in H\"ormander \cite{ho1}; the following definition is from Guillop\'e and Zworski \cite{gz}.  
\begin{defn} 
The \em $0$-regularized integral \em $\rint f$ of a smooth function (or density) $f$ on $X$ is defined, if it exists, as the finite part 
$$\rint_X  f :=  \underset{{\epsilon=0}}{\mathrm{f.p.}} \int_{x(p) > \epsilon} f(p) \dvol_g(p),$$
where $x$ is a boundary defining function.  
\end{defn}
The 0-regularized integral also appears in the literature as the 0-renormalized integral.  Guillop\'e and Zworski also defined the \em $0$-trace, \em which they used to formulate their Selberg trace formula for Riemann surfaces \cite{gz}.  

\begin{defn}
For an operator $A$ with smooth Schwartz kernel $A(z, y)$ on $X \times X$, the \em 0-trace \em of $A$ is 
$$\textrm{0-tr } (A) :=  \, \, \rint_X A(z, z) \dvol_g(z).$$
\end{defn} 
Both the $0$-regularized integral and the $0$-trace depend \`a priori on the choice of boundary defining function $x$ in their definitions.  However, in some cases, the 0-regularized integral is  independent of the boundary defining function.  The $0$-volume is defined to be the 0-regularized integral of the constant function $1$, 
$$\textrm{0-vol } (X) := \rint_X 1 \, \dvol_g.$$
Graham \cite{gr} showed that in even dimensions (by our convention, $n$ is odd), the $0$-volume is well defined and independent of the choice of boundary defining function used to define the $0$-regularized integral.  The 0-volume, also known as \em renormalized volume, \em of conformally compact manifolds is of general interest in conformal geometry; see for example Chang, Qing, and Yang \cite{cqy} and Chang, Gursky, and Yang \cite{cgy}.  

\subsection*{Bibliographical notes} 
The asymptotically cylindrical metrics or b-metrics were introduced by Melrose; an excellent reference for the geometric structure is \cite{tapsit}.  Asymptotically conic metrics also appear in the literature as scattering metrics, and asymptotically hyperbolic geometry has also been referred to as 0-geometry.  For each of these geometries, the resolvent operator has been constructed as an element of a pseudodifferential operator calculus.  These constructions  generalize H\"ormander's parametrix construction on closed manifolds to Riemannian manifolds with boundary whose metrics admit certain asymptotic forms near the boundary.  A novel feature of these constructions is the introduction of an extra symbol arising from the boundary.  The calculi for the asymptotically cylindrical and asymptotically hyperbolic geometries are known, respectively, as the b-calculus and the 0-calculus.  The 0-calculus was introduced by Mazzeo and Melrose in \cite{mmah}.  

Given an asymptotically hyperbolic metric, the existence of a boundary defining function $x$ such that the metric near the boundary has the form 
$$\frac{dx^2 + h(x, y, dy)}{x^2},$$
was first shown by Joshi and S\`a Barreto while studying inverse scattering results for asymptotically hyperbolic manifolds \cite{jsb}.  However, a simpler proof of this result may be found in \S 2 of Mazzeo and Pacard's work concerning Maskit combinations of Poincar\'e Einstein metrics \cite{mcpe}. 

M\"uller has also used the term ``asymptotically hyperbolic'' in his work concerning the spectral theory, geometry, and scattering theory of finite volume surfaces which have hyperbolic cusp type ends \cite{mu}.  Since those surfaces have finite volume, one may distinguish the asymptotically hyperbolic manifolds considered here if one notes that these are the asymptotically hyperbolic manifolds which are also \em conformally compact.  \em  

\section{Spectral and scattering theories}  

In 1984, Lax and Phillips \cite{lap} published an elegant and thorough study of the Laplace operator and its resolvent on convex cocompact hyperbolic manifolds, and in 1987, Perry \cite{pspec} generalized their results to the Laplacian with a short range potential.  These works provided useful insights to Mazzeo and Melrose \cite{mmah}, who studied the spectral theory of asymptotically hyperbolic manifolds.  An important contribution to their work was made by Guillarmou \cite{g0} nearly twenty years later.  We recall their main results below.  

\subsection{The resolvent}
Mazzeo and Melrose proved that the Laplacian on an asymptotically hyperbolic manifold of dimension $n+1$ has absolutely continuous spectrum, $\sigma_{ac} (\Delta) = [ \frac{n^2}{4}, \infty )$, and a finite pure point spectrum which is either empty or consists of finitely many eigenvalues, $\sigma_{pp} (\Delta) = \{\Lambda_k\}_{k=1} ^{m} \subset (0, \frac{n^2}{4} )$.  Mazzeo and Melrose constructed the resolvent kernel as an element of the so-called ``0-calculus,'' and showed that it is 
meromorphic on the half space
$$\left\{ \Re(s) > \frac{n}{2} \right\},$$
and admits a well-defined meromorphic extension to $\C$ with discrete poles of finite rank.  Above, the spectral parameter $s$ is related to the actual spectral parameter $\Lambda$ by 
$$\Lambda = s(n-s).$$  
The resonance set $\cR$ consists of the poles of the meromorphically continued resolvent 
$$R(s) := (\Delta - s(n-s))^{-1},$$ 
counted according to the multiplicity, 
$$\zeta \in \cR, \quad m (\zeta) :=   \textrm{ rank Res}_{\zeta} (\Delta - s(n-s))^{-1}.$$
This is not quite the full story; based on Borthwick and Perry's work concerning the scattering poles of asymptotically hyperbolic manifolds \cite{bp}, Guillarmou \cite{g0} investigated the behavior of the resolvent at the points $\left( \frac{n-k}{2} \right)_{k \in \N}$.  He made the important discovery that the presence of infinite rank poles at these points is determined by, and conversely determines, the type of expansion the metric admits near the boundary; this is all beautifully explained in \cite{g0}.  In fact, the resolvent admits a finite-meromorphic extension to $\C$ if and only if the metric $g$ admits an expansion of the form 
$$g = x^{-2} \left(dx^2 + \sum_{k = 0} ^{\infty} x^{2k} h_k (y; dy) \right),$$
near the boundary, where $x$ is a boundary defining function.  This is the case for many important examples, including convex cocompact hyperbolic manifolds and conformally compact manifolds which are hyperbolic near infinity.  In the latter case, the result is due to Guillop\'e and Zworski \cite{gzpol}. 

\subsection{The scattering operator}
The resolvent operator is intimately related to the scattering operator, whose definition we recall below.  
\begin{defn} \label{scat} 
Let $(X, g)$ be an asymptotically hyperbolic $n+1$ dimensional manifold with boundary defining function $x.$  For Re$(s) = \frac{n}{2}, \, s \neq \frac{n}{2},$ a function $f_1 \in \cC^{\infty} (\partial X)$ determines a unique solution $u$ of 
$$(\Delta - s(n-s) ) u = 0, \quad u \sim x^{n-s} f_1 + x^s f_2, \textrm{ as } x \to 0,$$
where $f_2 \in \cC^{\infty} ( \partial X).$  This defines the map called the \em scattering operator, \em  $S(s) : f_1 \mapsto f_2.$ 
\end{defn}

Heuristically, the scattering operator which is classically a scattering matrix, acts as a Dirichlet to Neumann map, and physically it describes the scattering behavior of particles.  Joshi and S\`a Barreto \cite{jsb} proved that the scattering operator extends meromorphically to $s \in \C$ as a family of pseudodifferential operators of order $2s - n$.  Renormalizing the scattering operator as follows gives a meromorphic family of Fredholm operators with poles of finite rank,   
$$\tilde{S} (s) : = \frac{\Gamma( s - \frac{n}{2})}{\Gamma(\frac{n}{2} - s) } \Lambda^{n/2 - s} S(s) \Lambda^{n/2 - s},$$
where 
$$\Lambda : = \frac{1}{2} (\Delta_h + 1)^{1/2}.$$
Above, $\Delta_h$ is the Laplacian on $\partial X$ for the metric $h(x) \big|_{x=0}.$  Note that this definition depends on the choice of boundary defining function $x$.  The multiplicity of a pole or zero of $S(s)$ is defined to be
$$\nu (\zeta) = - \tr [ \textrm{Res}_{\zeta} \tilde{S}'(s) \tilde{S}(s)^{-1} ].$$
The scattering multiplicities are related to the resonance multiplicities in \cite{bp} and \cite{g0} by 
$$\nu(s) = m(s) - m(n-s) + \sum_{k =1} ^{\infty}  \left( \chi_{n/2 - k} (s) - \chi_{n/2 + k} (s) \right) d_k,$$
where 
$$d_k = \dim \ker \tilde{S} (\frac{n}{2} + k),$$
and $\chi_{p}$ denotes the characteristic function of the set $\{p\}.$  This relationship is originally due to Borthwick and Perry \cite{bp}, but it was Guillarmou \cite{g0} who explicitly identified the numbers $d_k$ as the dimensions of the kernels of certain natural conformal operators acting on the conformal infinity.  These operators are related to the so-called ``GJMS'' operators \cite{gjms} and Q-curvature \cite{fg2}.  

\subsection{Counting estimates}
On a compact Riemannian manifold, the number of eigenvalues of the Laplacian grows according to Weyl's formula.  For convex cocompact hyperbolic manifolds of even dimension, Guillarmou \cite{gu1} obtained a Weyl-type asymptotic formula which interestingly depends on the topological entropy of the geodesic flow.  If $X$ is conformally compact and hyperbolic near infinity, upper and lower bounds for the resonance counting function were achieved through the work of Guillop\'e and Zworski \cite{gzpol} and completed by Borthwick \cite{bor} who demonstrated that the number of resonances (counted according to the spectral parameter $s$) in a ball of radius $r$ is $O(r^{n+1})$.  Moreover, the scattering resonance set satisfies the lower bound, 
$$ N^{\text{sc}}(r) \geq c B(X, g) r^{n+1} $$ 
where $c$ is a (universal) positive constant, and $B(X, g)$ is the 0-volume of $X$ if the dimension of $X$ is even or the Euler characteristic of $X$ if the dimension of $X$ is odd.  Counting estimates were  obtained for surfaces by Guillop\'e and Zworski \cite{gz2}.  Counting estimates for the scattering and resolvent resonances of asymptotically hyperbolic manifolds in general remains an intriguing and challenging open problem.  

\subsection*{Bibliographical notes} %%%
Guillarmou's Weyl law is in fact not the main result of \cite{g1}.  We highly recommend \cite{g1} to interested readers.  Guillarmou develops a Birman-Krein theory for asymptotically hyperbolic manifolds whose metrics have a certain form at the boundary.  These metrics include the natural examples:  convex cocompact hyperbolic manifolds and conformally compact manifolds hyperbolic near infinity.  He first shows that one may define the spectral measure via a 0-regularized integral and demonstrates its regularity properties on $\R$ and $\C$.  He goes on to study the determinant of the scattering operator and its relationship to the dynamical zeta function.  This work both provides new tools and suggests open problems on asymptotically hyperbolic manifolds.  For example, it would be interesting to generalize the Weyl law to these asymptotically hyperbolic manifolds.  

In \cite{lee} and \cite{lee2}, Lee made further contributions to the scattering theory and Fredholm properties of the Laplacian on  asymptotically hyperbolic manifolds; see also S\`a Barreto \cite{sb}.  Graham and Zworski \cite{grz} studied the scattering matrix in conformal geometry; this work is related to the numbers $d_k$ and the conformal infinity.  There is little known about the spectral theory of conformally compact manifolds which are not asymptotically hyperbolic: those which have variable curvature at infinity.  The only reference of which we are aware is Borthwick \cite{bor2}, who developed the scattering theory by treating the Laplacian as a degenerate elliptic operator with non-constant indicial roots.  The variability of the roots poses a significant challenge to the analysis because it implies that the resolvent admits only a partial meromorphic continuation.  

%%%%%%%%%%%%%%%%%%%%%%%%%%%%%%%%%%%%%%%%%%%%
\section{Dynamics} 
Classical mechanics are mathematically described by the dynamics of the geodesic flow.  We begin by recalling the definitions.  Let $SX$ denote the unit tangent bundle, and let $G^t$ be the geodesic flow on $SX$ at time $t$.  
\begin{defn} The \em geodesic flow \em on a Riemannian manifold $(X, g)$ is the map $G: SX \times \R \to SX$ such that for $v \in S_x X$, 
$$G^t (v) = \dot{\gamma_v} (t),$$
where $\gamma_v$ is the geodesic from $x \in X$ with initial tangent vector $v$, and $\dot{\gamma_v} (t)$ is its tangent vector at the point $\gamma_v (t)$.  The \em orbit \em of a vector $v \in S_x X$ for the geodesic flow is 
$$\{G^t (v) | t \in \R \}.$$
 \end{defn} 

The \em length spectrum \em is the set of lengths of closed geodesics.  We shall use 
$$\cL \quad \textrm{and} \quad \cL_p$$
to denote, respectively, the set of closed geodesics and the set of primitive closed geodesics.  Letting $\pi : SX \to X$ be the canonical projection, for all unit vectors $v \in SX$, the restriction of $\pi$ to the orbit of $v$ under $G$ is bijective onto the geodesic $\gamma_v \subset X$.  Thus, we see that the closed orbits of the flow $G$ are in bijection with the closed geodesics of $X$.  

Under certain conditions, the geodesic flow has desirable properties.  It is \em hyperbolic \em (see Anosov \cite{anosov}) if for each $\xi \in SX,$ $T(SX)_{\xi}$ splits into a direct sum 
\begin{equation} \label{split} T(SX)_{\xi} = E^s _{\xi} \oplus E^u _{\xi} \oplus E_{\xi},\end{equation}
where $E^s _{\xi}$ is exponentially contracting, $E^u _{\xi}$ is exponentially expanding, and $E_{\xi}$ is the one dimensional subspace tangent to the flow.  This definition was first given for closed manifolds as in Klingenberg \cite{kl}, but the same definition holds for open manifolds; see for example Bolton \cite{bolt}.  A closed $G^t$ invariant set $\Omega \subset X$ without fixed points is \em hyperbolic \em if the tangent bundle restricted to $\Omega$ is a Whitney sum,
$$T_{\Omega} X = E + E^s + E^u$$
of three $G^t$ invariant sub-bundles, where $E$ is the one dimensional bundle tangent to the flow, and $E^s$, $E^u$ are, respectively, exponentially contracting and expanding.  The \em Sinai-Ruelle-Bowen potential \em is a H\"older continuous function defined for $\xi \in SX$ by, 
\begin{equation} \label{srb} H(\xi) : = \frac{d}{dt} \big|_{t=0} \log \det dG^t |_{E^u _{\xi}}.\end{equation} 
This potential gives the instantaneous rate of expansion at $\xi$.  

A  vector $\xi \in SX$ is \em non-wandering \em for the geodesic flow if for all neighborhoods $U \subset SX$ with $\xi \in U$, there exists a sequence $\{ t_n\} \to \infty$ such that for all $n \in \N$,
$$U \cup G_{t_n} U \neq \emptyset.$$
These vectors form the \em non-wandering set.  \em  The geodesic flow satisfies Smale's \em Axiom A \em \cite{sm} if the non-wandering set $\Omega$ is hyperbolic, and the periodic points of the flow are dense in $\Omega$.  The flow is \em topologically transitive \em on $\Omega$, if for any open $U$ and $V \subset \Omega$, there exists $n > 0$ such that $G^n (U) \cap V \neq \emptyset$.  The definition of \em basic set \em is due to Bowen \cite{b74}.  A \em basic set \em $\Omega$ of an Axiom A flow $G^t$ on a Riemannian manifold $(X, g)$ is a hyperbolic set in which periodic orbits are dense, $G^t | \Omega$ is topologically transitive, and $\Omega = \cap _{t \in \R} G^t U$ for some open neighborhood $U \supset \Omega$.  An excellent reference for the general theory of dynamics is the book by Katok and Hassellblatt \cite{kh}.  

On a closed manifold with hyperbolic geodesic flow, the \em topological entropy of the geodesic flow \em is given by the limit
\begin{equation} \label{eq:ent1} h = \lim_{T \to \infty} \frac{\log \# \{ \gamma \in \cL : l(\gamma) \leq T \} }{T}, \end{equation} 
which was shown to exist in this setting by Bowen \cite{b72}.  When $X$ is compact and has no conjugate points, Freire and Ma\~n\'e \cite{frmane} showed that the topological entropy $h$ of the geodesic flow is equivalently given by the \em volume entropy, \em the log-volume growth rate in the universal covering, 
$$h = \lambda(X) := \limsup_{r \to \infty} \frac{ \log \Vol (B_r (x))}{r}.$$
Above, $\Vol(B_r(x))$ denotes the volume of the ball of radius $r$ and center $x$ in the universal covering of $X$; see also Manning \cite{man1}. 

For asymptotically hyperbolic manifolds, Joshi and S\`a Barreto showed that all closed geodesics are contained in a compact subset \cite{jsb2}.  Therefore, we may define the topological entropy of the geodesic flow by restricting to the non-wandering set.  First, we recall that for large $T$ and small $\delta > 0,$ a finite set $Y \subset SX$ is $(T, \delta)$ separated if, given $\xi, \xi' \in Y, \xi \neq \xi',$ there is $t \in [0, T]$ with $d(G^t \xi, G^t \xi') \geq \delta.$  Here the distance on $SX$ is given by the Sasaki metric.  We may then define the topological entropy of the geodesic flow as follows.   
\begin{defn}  Let $(X, g)$ be an asymptotically hyperbolic manifold with negative sectional curvatures.  Let $\Omega \subset SX$ be the non-wandering set of the geodesic flow.  We define the \em topological entropy of the geodesic flow \em to be the limit 
$$h :=  \lim_{\delta \to 0} \limsup_{T \to \infty} \frac{ \log \sup \# \{ Y \subset \Omega: Y \textrm{ is $(T, \delta)$ separated} \}}{T}.$$
\end{defn}  
In \cite{row1}, we demonstrate that this definition of topological entropy coincides with the definition given for the topological entropy of a hyperbolic geodesic flow on a closed manifold (\ref{eq:ent1}).  

For convex cocompact hyperbolic manifolds $\bbH^{n+1} / \Gamma$, Patterson \cite{pat} demonstrated that $h=\delta$ is the exponent of convergence for the Poincar\'e series for $\Gamma$, and Sullivan \cite{su} showed that $h$ is also dimension of the limit set of $\Gamma$.  Yue \cite{yue} has extended their results to some complete spaces of non-constant curvature.  These are known as ``convex cocompact manifolds'' and arise as the quotient of a complete manifold with pinched negative sectional curvatures by a convex cocompact group.  

The pressure of a function is a concept in dynamical systems arising from statistical mechanics which measures the growth rate of the number of separated orbits weighted according to the values of $f$; see for example Walters \cite{walt} and Manning \cite{man1}.  The rather cumbersome definition of pressure is nonetheless useful to generalize entropy, 
$$ \mathfrak{p}(f) = \lim_{\delta \to 0} \limsup_{T \to \infty}  \frac{\log \sup \left\{ \sum_{\xi \in Y} \exp \int_0 ^T f(G^t \xi) dt;  Y \textrm{ is $(T, \delta)$ separated} \right\}}{T}. $$
In the compact setting, the topological pressure of a function $f: SX \to \R$ 
$$ \mathfrak{p}(f) = \sup_{\mu} \left( h_{\mu} + \int f d\mu \right), $$
where the supremum is taken over all $G^t$ invariant measures $\mu,$ and $h_{\mu}$ denotes the \em measure theoretical entropy \em of the geodesic flow.  

The dynamical theory of open manifolds was initiated in the 1970s by several authors.  We note the work of Eberlein \cite{e}, \cite{e1}, \cite{e2}; Eberlein and O'Neill \cite{eo}; and Bishop and O'Neill \cite{bo}.  A further reference is Eberlein, Hamenst\"adt, and Schroeder \cite{ehs}.  The focus of those works is dynamics on \em visibility manifolds:  \em complete, open manifolds with non-positive sectional curvatures.  Important contributions to the dynamical theory of hyperbolic manifolds were made by Sullivan \cite{su}, Lalley \cite{lal}, Guillop\'e \cite{gu}, Patterson \cite{pat}, and Perry \cite{plength}.  

\subsection{Dynamics of asymptotically hyperbolic manifolds} 
The following dynamical lemma was proven in \cite{row1} using methods suggested by Manning \cite{mann}, and we expect that it is known to experts, if not already in the literature.  Nonetheless, since it is a useful tool for dynamics on any complete manifold with pinched negative sectional curvatures, we include its statement here.    
\begin{lemma}  \label{le:lyapunov} Let $(X, g)$ be a smooth, complete, $n+1$ dimensional Riemannian manifold whose sectional curvatures $\kappa$ satisfy 
$$-k_1^2 \leq \kappa \leq - k_2 ^2$$
for some $0 < k_2 \leq k_1.$  Then, the Poincar\'e map about a closed orbit $\gamma$ of the geodesic flow has $n$ expanding eigenvalues $\{\lambda_i\}_{i=1} ^n$ which satisfy 
$$e^{k_2 l(\gamma)} \leq |\lambda_i| \leq  e^{k_1 l(\gamma)}  \textrm{ for } i=1, \ldots , n,
$$ and $n$ contracting eigenvalues $\{\lambda_i\}_{i=n+1} ^{2n}$ which satisfy  
$$e^{- k_1 l(\gamma)} \leq |\lambda_i | \leq e^{-k_2 l(\gamma)} \textrm{ for } i = n+1 , \ldots, 2n,$$
where $l(\gamma)$ is the period (or length) of $\gamma.$  
\end{lemma} 
The proof uses the Rauch Comparison Theorem, which can be found in do Carmo \cite{doc}, and basic properties of the Lyapunov exponents found in Luis-Pesin \cite{lup}.  

In \cite{row1}, we demonstrated a ``prime orbit theorem'' for the geodesic flow of negatively curved asymptotically hyperbolic manifolds.  

\begin{theorem} \label{theorem3} Let $(X, g)$ be an asymptotically hyperbolic $n+1$ dimensional manifold with negative sectional curvatures.  If the the topological entropy of the geodesic flow $h>0$, then the length spectrum counting function 
\begin{equation} \label{eq:lect} N(T) := \# \{ \gamma \in \cL : l(\gamma) \leq T \} \textrm{ satisfies } \lim_{T \to \infty} \frac{T N(T)}{e^{hT}} = 1. \end{equation} 
\end{theorem} 

The proof of this result relies on the following key ingredients.  First, Joshi and S\`a Barreto demonstrated that the closed geodesics of \em any \em asymptotically hyperbolic manifold are contained in a compact subset.  Next, since the sectional curvatures are negative, it is straightforward to generalize the ``separation lemma'' of Jakobson, Polterovich, and Toth \cite{jpt} to asymptotically hyperbolic manifolds; see Lemma 3.4 of \cite{row1}.  We use this lemma and the dynamical lemma above in combination with the aforementioned references on the dynamics of visibility manifolds to show that, if the topological entropy is positive, then the non-wandering set for the geodesic flow is a basic set.  It follows from Eberlein \cite{e1} that the geodesic flow restricted to the non-wandering set is an Axiom A flow restricted to a basic set, so we apply the main result of Parry and Pollicott \cite{pot} to the dynamical zeta function which completes the proof of the theorem.  In fact, we are currently exploring further implications for the spectral theory of asymptotically hyperbolic manifolds with negative sectional curvatures based on more refined properties of the dynamical zeta funtions.  

\subsection{Dynamical zeta functions}
There are two key objects which relate the spectral and dynamical theories.  The first is the dynamical zeta function, which is to the geodesic length spectrum as the Riemann zeta function is to the prime numbers.  Let 
\begin{equation} \label{dz} Z(s) = \exp  \left( \sum_{\gamma \in \cL_p} \sum_{k=1} ^{\infty} \frac{e^{-k s l_p(\gamma)}}{k} \right), \end{equation}
where $\cL_p$ consists of primitive closed orbits of the geodesic flow and $l_p(\gamma)$ is the primitive period (or length) of $\gamma \in \cL_p$.  Patterson and Perry \cite{pp} defined the following weighted dynamical zeta function, 
$$ \tilde{Z} (s) = \exp \left( \sum_{\gamma \in \cL_p} \sum_{k =1} ^{\infty} \frac{e^{-k s l_p(\gamma)}}{k \sqrt{| \det(I - \cP^k _{\gamma}) |}} \right),$$
where $\cP^k _{\gamma}$ is the $k$-times Poincar\'e map of the geodesic flow around the (primitive) closed orbit $\gamma$.  The weighted dynamical zeta function is particularly interesting for its connections to the resonances of the resolvent; Perry \cite{p} and Guillarmou-Naud \cite{gn}, via the Hadamard factorization of the weighted dynamical zeta function, proved the Selberg trace formulae for convex cocompact hyperbolic manifolds.   Borthwick, Judge, and Perry \cite{bjp} used Birman-Krein theory and explicit model operators to describe the poles and zeros of the dynamical zeta function on hyperbolic surfaces with finite geometry.  

For asymptotically hyperbolic manifolds with negative sectional curvatures, we obtained the following result in \cite{row1} for the dynamical zeta functions.  
\begin{theorem} 
Let $(X, g)$ be an asymptotically hyperbolic $n+1$ dimensional manifold with negative sectional curvatures.  Then, the dynamical zeta function $Z(s)$ converges absolutely if and only if Re$(s) > h$, where $h$ is the topological entropy of the geodesic flow.  Moreover, if $h > 0$, then $Z$ admits a nowhere vanishing analytic extension to an open neighborhood of Re$(s) \geq h$ except for a simple pole at $s=h$.  The weighted dynamical zeta function, $\tilde{Z} (s)$ converges absolutely if and only if Re$(s) > \fp(-H/2)$, where $\fp$ is the topological pressure, and $H$ is the Sinai-Ruelle-Bowen potential. \end{theorem}

Our proof was facilitated by Joshi and S\`a Barreto's results which allow one to apply local estimates from dynamics on closed manifolds, and in particular Chen and Manning \cite{cm}, Bowen \cite{bo}, and Franco \cite{fr}.  

\subsection*{Bibliographical notes} 
Anosov's work was motivated by S.~Smale's lectures \cite{sm}; Anosov attended the lectures and subsequently answered Smale's conjectures in an impressive doctoral thesis \cite{anosov}.  An Anosov flow is a special example of Smale's Axiom A flows; Anosov's original definition was a flow that satisfies  ``condition U'' \cite{anosov}.  Another related flow is the Bernoulli flow which was studied by Lalley, who in that setting proved the prime orbit theorem \cite{lal}.  The prime orbit theorem for surfaces with totally geodesic boundary was proven by Guillop\'e \cite{gu}.  For convex cocompact hyperbolic manifolds in all dimensions, the prime orbit theorem is due to Perry \cite{plength}.  Guillarmou and Naud \cite{gn} obtained an estimate for the remainder term in Perry's prime orbit theorem.  

Patterson and Perry studied the dynamical zeta function in even dimensions \cite{pp}.  Their work includes an appendix by C. Epstein which demonstrates an interesting geometric property for conformally compact manifolds with constant negative curvature near infinity.  The metric admits a particular expansion in this neighborhood of infinity, and one may compute the coefficients in terms of the Riemannian curvature tensor and its derivatives.  
%%%%%%%%%%%%%%%%%%%%%%%%%%%%%%%%%%%%%
\section{Quantum and classical mechanics}
The movement of very small particles is described by \em quantum mechanics, \em the name of which derives from the observation that certain physical quantities such as the energy of an electron bound into an atom or molecule can only assume discrete values or \em quanta.  \em  One of the first mathematical formulations of quantum mechanics is matrix mechanics, due to Heisenberg, Born, and Jordan \cite{bjh} in 1925.  This theory was shortly followed by wave mechanics, introduced by Schr\"odinger \cite{swave} in 1926.  Two years later, Dirac \cite{dir} introduced transformation theory to unify and generalize these formulations and to unify the particle-wave duality of energy and matter observed in photons and electrons.  In these theories, eigenvalues of certain operators are used to describe the energy states of quantum particles.  One of the fundamental operators considered is the Laplace operator on a Riemannian manifold.  When the manifold has infinite volume, the \em resonances \em of the resolvent describe the quantum states.  

The relationship between quantum and classical mechanics remains in many ways an open question.  One theory is that all objects obey the laws of quantum mechanics, and that classical mechanics is simply the quantum mechanics of a very large collection of particles.  Thus, the laws of classical mechanics ought to follow from the laws of quantum mechanics and limits of large quantum systems.  However, this breaks down with certain chaotic systems; see for example  Gutzwiller \cite{gut}.  A further crux is the Einstein-Podolsky-Rosen paradox \cite{epr}; simply put, the laws of quantum mechanics would seem to violate the most natural and basic laws of classical mechanics.  Mathematical physicists would like to reconcile quantum and classical mechanics, so it is helpful to understand interactions between them.  The wave group is one bridge between the quantum and classical mechanics.  

\subsection{The wave group}
The (even) wave kernel is the Schwartz kernel of the fundamental solution to 
$$\left( \partial_t ^2 + \Delta - \frac{n^2}{4} \right)U(t, w, w') = 0,$$
$$U(0, w, w') = \delta(w-w'), \quad \frac{\partial}{\partial t}  U(0, w, w') = 0.$$
Due to the semi-group property with respect to time, the wave kernel is also referred to as the wave group with the notation 
$$\cos \left( t \sqrt{ \Delta - n^2/4 }\right).$$
The wave group on asymptotically hyperbolic manifolds was constructed by Joshi and S\`a Barreto in  \cite{jsb2} as an element of an operator calculus defined on a certain manifold with corners obtained by blowing up $\R^+ \times X \times X$ along two submanifolds (the diagonal is blown up first, followed by the submanifold where the diagonal meets the corner).    

They demonstrated that the wave group for an asymptotically hyperbolic manifold has a well-defined  0-trace, and that the singular support of 0-tr $\cos (t \sqrt{ \Delta - n^2 /4} )$ is contained in the set of lengths of closed geodesics.  In effect, they generalized Duistermaat and Guillemin \cite{dg} to the asymptotically hyperbolic setting.  

\subsection{Trace formulae}
The leading term in the following trace formula follows immediately from Joshi and S\`a Barreto \cite{jsb2}.  The main point of the result is the long time remainder estimate which allows one to use the trace formula to uncover a rapport between the size of the topological entropy and the existence of pure point spectrum thereby establishing a quantitative relationship between the classical and quantum mechanics.  This long-time estimate is based on the techniques of Jakobson, Polterovich, and Toth \cite{jpt} using the stationary phase method augmented by estimates using the iterative construction of the wave group in \cite{jsb2} and properties of the wave group demonstrated by B\'erard \cite{ber}.  The estimate is in \em Ehrenfest time, \em in which the oscillation of the test function and the size of its support satisfy a certain relationship; see for example Zelditch \cite{zel}.  

\begin{theorem}  \label{theorem1} 
Let $(X, g)$ be an asymptotically hyperbolic $n+1$ dimensional manifold with negative sectional curvatures.  As a distributional equality in $\cD'((0, \infty))$, 
\begin{equation} \label{trace1} \textrm{0-tr} \cos ( t \sqrt{ \Delta - n^2/4}) = \sum_{\gamma \in \cL_p} \sum_{ k =1} ^{\infty} \frac{ l(\gamma) \delta(t - k l(\gamma))}{ \sqrt{ | \det(I - \cP^k _{\gamma}) | }} + A(t).\end{equation}
The remainder $A$ is smooth and satisfies the following estimate:  there exists $\epsilon > 0$ such that for any $t_0 > 0$, there exists a $C>0$ such that 
\begin{equation}\label{eq:longt} \left| \int_0 ^{\infty} A(t) \cos (\lambda t) \rho(t) dt \right| \leq C \end{equation}
for all $\lambda >1$ and $\rho \in \cC^{\infty} _0 ([t_0, \epsilon \ln \lambda])$.  The constant $C$ depends only on $t_0$ and $|| \rho ||_{\infty}$; $C$ is independent of $\lambda$.  
 \end{theorem}
 
On a convex cocompact hyperbolic manifold $X^{n+1} = \bbH^{n+1}/ \Gamma$, Perry \cite{p} and Guillarmou Naud \cite{gn} proved the Selberg trace formula.  The dynamical side of their trace formula is the following distributional equality in $\cD'((0, \infty))$,
\begin{equation} \label{eq:cctrace} \textrm{0-tr} \cos ( t \sqrt{ \Delta - n^2/4}) = \sum_{\gamma \in \cL_p} \sum_{ k =1} ^{\infty} \frac{ l(\gamma) \delta(|t| - k l(\gamma))}{ \sqrt{ | \det(I - \cP^k _{\gamma}) | }} + A(X) \frac{\cosh{\frac{t}{2}}}{(\sinh{\frac{|t|}{2}})^{n+1}}, \end{equation}
where
$$A(X) = \left\{ \begin{array}{ll} n!! 2^{-\frac{3(n+1)}{2}} (-\pi)^{-\frac{n+1}{2}} \textrm{0-vol}(X) & \textrm{ if $n+1$ is even} \\ 0 & \textrm{if $n+1$ is odd.} \end{array} \right.$$
In the case of surfaces, the trace formula was originally demonstrated by Guillop\'e and Zworski \cite{gz}; that result is somewhat more general, since they allow the surface to have cusps.  

It is interesting to note that the remainder in this Selberg trace formula has markedly different properties from the remainder in the Duistermaat Guillemin trace formula for closed manifolds \cite{dg}.  Specifically, there are examples of convex co-compact hyperbolic manifolds whose geodesic flow has topological entropy $h > \frac{n}{2}$; see for example Canary, Minsky, and Taylor \cite{cmt} who construct examples of hyperbolic 3-manifolds.  By the prime orbit theorem for these manifolds \cite{plength}, Lemma \ref{le:lyapunov}, and a straightforward estimate, the renormalized wave trace lies only in $\cD'((0, \infty))$ rather than $\cS '((0, \infty))$.  On the other hand, the wave trace on a compact manifold is always a Schwartz distribution.  

Borthwick's Poisson formula \cite{bor2} was the key to our resonance wave trace formula \cite{row1}.  
 \begin{corollary} \label{corollary1} 
Let $(X, g)$ be an $n+1$ dimension conformally compact manifold hyperbolic near infinity with negative sectional curvatures.  Then, we have the distributional equality in $\cD'((0, \infty))$, 
\begin{equation} \label{trace2} \frac{1}{2} \sum_{s \in \cR } m(s) e^{(s - n/2)t} + \frac{1}{2} \sum_{k =1} ^{\infty} d_k e^{-kt} = \sum_{\gamma \in \cL_p} \sum_{ k =1} ^{\infty} \frac{ l(\gamma) \delta(t - k l(\gamma))}{ \sqrt{ | \det(I - \cP^k _{\gamma})| }} + C(t). \end{equation}
Above, the numbers $d_k$ are determined by natural conformal operators acting on the conformal infinity; see \cite{gn}.  The long time asymptotics of the remainder are given by (\ref{eq:longt}).
\end{corollary}

\subsection{Interactions between quantum and classical mechanics}

Our final and most interesting result in \cite{row1} is a quantitative relationship between the topological entropy of the geodesic flow and the existence of pure point spectrum, in the spirit of Patterson \cite{pat} and Phillips \cite{ph}.  The key ingredients in the proof are the long-time remainder estimate in the trace formula and a Littlewood counting estimate used in analytic number theory; see for example Rubinstein and Sarnak \cite{sarn}, Phillips and Rudnik \cite{phr}, and Karnauk \cite{karn}.  The idea was inspired by Jakobson, Polterovich, and Toth \cite{jpt} who used this counting technique (which they called the \em Dirichlet box principle\em) to estimate the remainder in Weyl's law on surfaces with pinched negative curvature.  

\begin{corollary} \label{corollary3}
Let $(X, g)$ be an $n+1$ dimension conformally compact manifold hyperbolic near infinity with negative sectional curvatures and topological entropy $h$ for the geodesic flow.  Let $0 < k_2 \leq 1 \leq k_1$ be such that the sectional curvatures $\kappa$ satisfy $-k_1 ^2 \leq \kappa \leq -k_2 ^2$.  If $h > \frac{nk_1}{2},$ then $\sigma_{pp} (\Delta) \neq \emptyset,$ and moreover, there is $\Lambda_0 = s_0 (n-s_0) \in \sigma_{pp} (\Delta)$ with $s_0 \geq h + n(1-k_1)/2.$  Conversely, if $h \leq \frac{nk_2}{2},$ then $\sigma_{pp}(\Delta) = \emptyset.$  
\end{corollary}
The contrapositive of the corollary implies 
$$\sigma_{pp}(\Delta) \neq \emptyset \Rightarrow h > \frac{nk_2}{2}, \qquad \sigma_{pp}(\Delta) = \emptyset \Rightarrow h \leq \frac{nk_1}{2}.$$ 

\subsection*{Bibliographical notes} 
The relationship between the pure point spectrum and the topological entropy of the geodesic flow in \cite{row1} is not sharp like the constant curvature case.  It would be interesting to construct examples and study the effects of varying the metric and the curvature.  In the case of conformally compact metrics hyperbolic near infinity, our result seems to imply that if the topological entropy is large, one cannot destroy pure point spectrum by perturbing the metric on a compact set, or that such a perturbation would necessarily change the curvature pinching constants.  What are the physical implications?  For readers interested in the physical aspects, we recommend Naud's work concerning Riemann surfaces \cite{n}.  

\section{The horizon} %%%%%%%%%%%%%%%
During the preparation of this survey, Borthwick and Perry applied our trace formulae \cite{row1} to obtain inverse scattering results for conformally compact manifolds hyperbolic near infinity \cite{bp2}.  Their strongest result is in the case of surfaces which have a common resonance set and the same model at infinity; the set of all such surfaces is compact in the $\cC^{\infty}$ topology.  In higher dimensions, they assume the sectional curvatures are negative and depending on whether the dimension is even or odd, they assume either a common resonance set or a common scattering phase and obtain compactness results.  It may be possible to obtain further results via our trace formulae and Corollary \ref{corollary3}; in any case, there are many open problems for the spectral theory and dynamics of asymptotically hyperbolic manifolds.  We conclude this survey with the following conjectures and ideas.  

\subsection{Wave trace remainder estimates}  

On a compact manifold $X$ with Laplace spectrum $\{ \lambda_k \}_{k=1} ^{\infty}$, the wave trace is \em formally, \em 
\begin{equation}  \sum_{k =1} ^{\infty} e^{i \sqrt{\lambda_k} t} = \sum_{\gamma \in \cL_p} \sum_{ k =1} ^{\infty} \frac{l(\gamma)\delta(|t| - k l(\gamma))}{\sqrt{|\det( I - \cP_{\gamma} ^k)|}} + A(t) \in \cS'((0, \infty)). \end{equation}
The principle term in the dynamical side,
\begin{equation} \label{eq:cont1} \sum_{\gamma \in \cL_p} \sum_{ k=1} ^{\infty} \frac{l(\gamma)\delta(|t| - k l(\gamma))}{\sqrt{|\det( I - \cP_{\gamma} ^k)|}} \end{equation}
may have exponential growth depending on the curvature bounds and the topological entropy of the geodesic flow.  Since its sum with the remainder term $A(t)$ lies in $\cS'((0, \infty))$, the remainder term cannot lie in $\cS'((0, \infty))$ unless the principle term (\ref{eq:cont1}) is also in $\cS'((0, \infty))$.  In contrast, the remainder in the renormalized wave trace for convex co-compact hyperbolic manifolds (\ref{eq:cctrace}) has exponential decay for large $t$.  It is therefore not unreasonable to expect the following.  

\begin{conj} \label{con1}
Let $(X, g)$ be an asymptotically hyperbolic $n+1$ dimensional manifold with negative sectional curvatures.  As a distributional equality in $\cD'((0, \infty))$, 
\begin{equation} \label{trace1} \textrm{0-tr} \cos ( t \sqrt{ \Delta - n^2/4}) = \sum_{\gamma \in \cL_p} \sum_{ k =1} ^{\infty} \frac{ l(\gamma) \delta(t - k l(\gamma))}{ \sqrt{ | \det(I - \cP^k _{\gamma}) | }} + A(t).\end{equation}
The remainder $A$ is smooth on $(0,\infty)$, and there exist $C$, $k > 0$ which depend only on $n$, $k_1$, and $k_2$ such that 
$$|A(t)| \leq Ce^{kt} \quad \forall \quad t > 1,$$
where the sectional curvatures $\kappa$ satisfy $-k_1 ^2 \leq \kappa \leq -k_2 ^2$.  
\end{conj}

We propose a new approach to wave trace remainder estimates using the higher wave invariants.  It is unclear how feasible this approach is, however, it is of independent interest that Zelditch's results \cite{z2} extend to asymptotically hyperbolic manifolds.  

On a compact manifold, the singular support of the wave trace for $t>0$ is contained in the set of lengths of closed geodesics.  The wave trace has the singularity expansion 
$$\textrm{tr} \cos(t \sqrt{\Delta}) = e_0 (t) + \sum_{\gamma} e_{\gamma} (t),$$
where
$$e_{\gamma} (t) \sim a_{-1}(\gamma) (t - l(\gamma) + i 0)^{-1} + \sum_{k=0} ^{\infty} a_k (\gamma) (t-l(\gamma) + i0)^k \log (t-l(\gamma) + i0).$$
The wave invariants are the coefficients $a_k (\gamma)$ in the singularity expansion of the wave trace.  By \cite{jsb2} and \cite{dg}, the renormalized wave trace for an asymptotically hyperbolic manifold has an analogous singularity expansion 
$$\textrm{0-tr} \cos ( t \sqrt{ \Delta - n^2/4}) \sim e_0 (t) + \sum_{\gamma} e_{\gamma} (t).$$
The principle wave invariant,
$$a_{-1}(\gamma) = \frac{ l_p(\gamma)} { \sqrt{ | \det(I - \cP _{\gamma}) | }},$$
where $l_p(\gamma)$ is the length of the primitive period of $\gamma$.  The higher wave invariants $a_{k}$ for $k \geq 0$ were studied by Zelditch \cite{z1}, \cite{z2}; and Guillemin \cite{gu1}, and \cite{gu2}.  Recall that a closed geodesic is \em non-degenerate \em if the Poincar\'e map $\cP_{\gamma}$ at $\gamma$ is any symplectic sum of non-degenerate elliptic, hyperbolic, or loxodromic parts.  The elliptic eigenvalues come in complex conjugate pairs of modulus one, $e^{\pm i \alpha_j}$; the hyperbolic eigenvalues come in inverse pairs of real eigenvalues, $e^{\pm \lambda_j}$ (these could also be negative); and the loxodromic eigenvalues come in 4-tuplets, $e^{\pm \mu_j \pm i\nu_j}$; where the \em Floquet exponents, \em $\alpha_j, \lambda_j, \mu_j, \nu_j \in \R$.  The non-degeneracy condition is equivalent to the independence over $\Q$ of the Floquet exponents together with $\pi$.  To study the wave invariants, it is useful to introduce Fermi normal coordinates $(s, y)$ along the geodesic $\gamma$.  Invariant contractions against $\frac{\pa}{\pa s}$ and against the Jacobi eigenfields $Y_j$, $\overline{Y}_j$, with coefficients given by invariant polynomials in the components $y_{jk}$ are \em Fermi-Jacobi polynomials.  \em  The \em Floquet invariants, \em $\beta_i = (1-\rho_i)^{-1}$, where $\{\rho_i\}_{i=1} ^{2n+2}$ (in dimension $n+1$) are the eigenvalues of $\cP_{\gamma}$.  

The first main result of \cite{z2} is that the wave invariants may be written in the form 
$$a_{\gamma k} = \mathcal{F}_{k, -1} (D) \cdot Ch(x) |_{x=\cP_{\gamma}},$$
where 
$$Ch(x) = \frac{i^{\sigma}}{\sqrt{|\det (I-x)|}}$$
is a character of the metaplectic representation (with $\sigma$ a certain Maslov index), and where $\mathcal{F}_{k, -1} (D)$ is an invariant partial differential operator on the metaplectic group $Mp (2(n+1), \R)$ which is canonically fashioned from the germ of the metric $g$ at $\gamma$.  In particular, upon close inspection of \cite{z2}, we have the following.  Note that the statement of the result 

\begin{theorem}  \label{th:wave1} Let $(X, g)$ be an asymptotically hyperbolic $n+1$ dimensional manifold, and let $\gamma$ be a non-degenerate closed geodesic.  The wave invariants, 
$$a_k (\gamma) =  \int_{\gamma} I_{\gamma; k} (s; g) ds,$$
where $I_{\gamma; k}$ satisfies 
\begin{enumerate}[i.]
\item $I_{\gamma; k} (s,g)$ is a homogeneous Fermi-Jacobi-Floquet polynomial of weight $-k-1$ in the data $\{y_{ij}, \dot{y}_{ij}, D^{\beta} _{s,y} g \}$ with $|\beta| \leq 2k +4$;
\item The degree of $I_{\gamma;k}$ in the Jacobi field components is at most $6k + 6$; 
\item At most $2k+1$ indefinite integrations over $\gamma$ occur in $I_{\gamma; k}$;
\item The degree of $I_{\gamma; k}$ in the Floquet invariants $\beta_j$ is at most $k+2$.  
\end{enumerate}
\end{theorem}

\begin{remark}  The exact expression for the wave invariants is given in \S 5 of \cite{z2}.  Estimating these expressions with respect to $k$ and applying Lemma \ref{le:lyapunov} is one approach to Conjecture \ref{con1}.  
\end{remark}

The wave invariants at a closed geodesic $\gamma$ are expressed by Guillemin  \cite{gu1}, \cite{gu2} and Zelditch \cite{z1}, \cite{z2} as non-commutative residues of the wave group and its time derivatives at $t= l(\gamma)$,
$$a_{\gamma k} = \res D^{k} _t e^{it\sqrt{\Delta}} |_{t=l(\gamma)} := \Res_{s=0} \Tr D^k _t e^{it\sqrt{\Delta}} \Delta ^{-\frac{s}{2}} |_{t=l(\gamma)}.$$
Since this residue is invariant under conjugation by microlocal unitary Fourier integral operators, the wave invariants may be calculated by putting the wave group into a microlocal quantum Birkhoff normal form around $\gamma$ and by determining the residues of the resulting wave group of the normal form.  Both Guillemin \cite{gu2} and Zelditch \cite{z1} computed the form for elliptic closed geodesics; subsequently, Zelditch \cite{z2} computed the form for non-degenerate closed geodesics.  In this case, the form includes transverse elliptic harmonic oscillators $\hat{I}^e _j = \frac{1}{2} (D_{y_j} ^2 + y_j ^2)$, real hyperbolic action operators $\hat{I}^h _j$, and loxodromic action operators $\hat{I}^{ch, Re} _j$, $\hat{I}^{ch, Im} _j$.  

\begin{theorem} \label{th:wave2} There exists a microlocally elliptic Fourier integral operator $W$ from the conic neighborhood of $\R^+ _{\gamma}$ in $T^* (N_{\gamma})$ to the corresponding cone in $T^* _+ \bbS^1$ in $T^* (\bbS^1 \times \R^{n+1})$ such that 
$$W \sqrt{\Delta} W^{-1} \equiv D_s + \frac{1}{l(\gamma)} \left[ \sum_{j=1} ^p \alpha_j \hat{I}_j ^e + \sum_{j=1} ^q \lambda_j \hat{I}_j ^h + \sum_{j=1} ^c \mu_j \hat{I}_j ^{ch, Re} + \nu_j \hat{I}_j ^{ch, Im} \right]$$
$$+ \frac{p_1 (\hat{I}_1 ^e, \ldots, \hat{I}_p ^e, \hat{I}_1 ^h , \ldots , \hat{I}^h _q, \hat{I}^{ch, Re} _1, \hat{I}^{ch, Im} _1, \ldots , \hat{I}^{ch, Re} _c, \hat{I}^{ch, Im} _c )}{D_s} + \ldots$$
$$+\frac{p_{k+1} (\hat{I}^e _1, \ldots, \hat{I}^{ch, Im} _c )}{D^k _s} + \ldots$$
where the numerators $p_l (\hat{I}^e _1, \ldots, \hat{I}^e _p, \hat{I}^h _1, \ldots , \hat{I}^{ch, Im} _c)$ are polynomials of degree $l+1$ in the variables $(\hat{I}^e _1, \ldots, \hat{I}^{ch, Im} _c)$, and the $k^{th}$ remainder term lies in the space $\oplus_{j=0} ^{k+2} O_{2(k+2-j)} \Psi^{1-j}$.  Here, $O_m \Psi^r$ is the space of pseudodifferential operators of order $r$ whose complete symbols vanish to order $m$ at $(y, \eta) = (0,0)$.  
\end{theorem}

Since the preceding results imply that the wave invariants at a non-degenerate closed geodesic determine, and conversely are determined by the Birkhoff normal form of the Laplacian, we have the following.  

\begin{theorem} \label{th:wave3} 
Let $(X, g)$ be an asymptotically hyperbolic $n+1$ dimensional manifold, and let $\gamma$ be a non-degenerate closed geodesic so that the only closed geodesics of the same length are $\gamma, \gamma^{-1}$.  Then the entire quantum Birkhoff normal form around $\gamma$ (and hence the classical Birkhoff normal form) is determined by the wave invariants of $\gamma$ and its iterates, and conversely, the wave invariants of $\gamma$ and its iterates are determined by the quantum normal form coefficients.
\end{theorem}

\textbf{Proofs:  }  
To prove Theorems \ref{th:wave1}--\ref{th:wave3}, it suffices to recall that by \cite{jsb2}, the closed geodesics of an asymptotically hyperbolic manifold lie in a compact subset. The proofs then follow immediately from the proofs of Theorem I, Theorem B, and Theorem II of \cite{z2}.  Rather than reproducing these arguments, we refer the reader to \cite{z2}.  In fact, if $(X, g)$ is any Riemannian manifold of dimension $n+1$, and $\gamma$ is a non-degenerate closed geodesic contained entirely in the interior of $X$, these results also hold. 

\qed 

\subsection{Spectral zeta functions}

The inverse spectral results in Guillemin \cite{gu1}, \cite{gu2}; and Zelditch \cite{z1}, and \cite{z2} use a characterization of the wave invariants via the spectral zeta function.  Recall that on a compact Riemannian manifold $(X, g)$ with non-zero Laplace spectrum $\{\lambda_k \}_{k =1} ^{\infty}$,
$$\zeta(s) := \sum_{k=1} ^{\infty} \lambda_k ^{-s}.$$
By the Weyl asymptotic formula, the zeta function is well defined and holomorphic when Re$(s) \gg 0$.  Based on the relationship between $\zeta$ and the heat kernel via the Mellin transform (see, for example \cite{rose}), it follows that $\zeta$ extends to a meromorphic function on $\C$.  Moreover, $s=0$ is a regular value which is used to define the \em zeta regularized determinant of the Laplacian, \em 
$$\det \Delta := e^{-\zeta ' (0) }.$$
Guillemin \cite{gu2} introduced the the zeta function 
$$\zeta_l (s) := \textrm{tr} ( e^{il \sqrt{\Delta}}) \Delta^{\frac{s}{2}},$$
which he used to characterize the wave invariants.  Guillemin's zeta function is well defined and holomorphic when Re$( -s) \gg 0$ and extends to a meromorphic function on $\C$ with poles at $\{-1, 0\} \cup \N$.  Zelditch observed that the wave trace invariants are the residues of this zeta function at its poles; this observation played a key role in the proofs of both Guillemin's \cite{gu2} and Zelditch's \cite{z1}, \cite{z2} inverse spectral results.  

On an asymptotically hyperbolic manifold, resonances play the role of the eigenvalues $\{\lambda_k \}$.  The ``scattering resonance set'' defined in \S 1 of Borthwick \cite{bor2} is, 
$$\cR_{sc}  := \{ s \in \cR \} \cup \bigcup_{k=1} ^{\infty} \left\{ \frac{n}{2} - k \textrm{ with multiplicity } d_k \right\}.$$
We define the ``scattering zeta function,'' 
$$\zeta_{sc} (s) := \sum_{\lambda \in \cR_{sc}, \, \lambda \neq n} \frac{1}{(n-\lambda)^s}.$$
A preliminary step in the study of the spectral zeta function and its relation to wave invariants for asymptotically hyperbolic manifolds is the following.  
\begin{lemma}  Let $(X, g)$ be a conformally compact manifold hyperbolic near infinity of dimension $n+1$.  Then, the   
scattering zeta function converges absolutely for Re$(s) > n+1.$
\end{lemma}
\textbf{Proof:  }  The proof follows immediately from Borthwick's counting estimate \cite{bor2}, 
$$\# \{ \lambda \in \cR_{sc} : | \lambda | \leq r \} = O(r^{n+1}).$$  
\qed

One may also use the $0$-regularized integral to define a ``wave zeta function'' in the spirit of Guillemin \cite{gu2}, 
$$\zeta_w (l,s) := 0\textrm{-tr} ( e^{il \sqrt{\Delta}}) \Delta^{\frac{s}{2}}.$$
In light of Guillarmou's Weyl law \cite{g1}, Borthwick's counting estimates and preliminary results for the ``relative zeta function'' \cite{bor2} which is related to M\"uller's spectral zeta function \cite{mu}, it is not unreasonable to expect results in the spirit of M\"uller generalize to conformally compact manifolds hyperbolic near infinity.  It may also be possible to extend the inverse spectral results for the wave invariants \cite{gu2} to these spaces via the wave zeta function.  These spectral zeta functions, the dynamical zeta functions, counting estimates for asymptotically hyperbolic manifolds, the Poisson formula for asymptotically hyperbolic manifolds, and the dynamical theory of conformally compact manifolds are some of many open problems to be explored.  

\section*{Acknowledgments}
I would like to thank G\'erard Besson and the organizers of the International Conference on Spectral Theory and Geometry, June 1--5, 2009 at the Institut Fourier.  I am also grateful to Yves Colin de Verdi\`ere, Vincent Grandjean, Werner M\"uller, Peter Sarnak, Samuel Tapie, and Steve Zelditch for discussions and correspondence, and to the anonymous referee for constructive comments.  

\nocite{*}
\bibliographystyle{cdraifplain}

\end{document}